\input amstex
\magnification 1200
\loadmsbm
\parindent 0 cm

\define\nl{\bigskip\item{}}
\define\snl{\smallskip\item{}}
\define\inspr #1{\parindent=20pt\bigskip\bf\item{#1}}
\define\iinspr #1{\parindent=27pt\bigskip\bf\item{#1}}
\define\einspr{\parindent=0cm\bigskip}

\define\co{\Delta}
\define\st{$^*$-}
\define\ot{\otimes}

\input amssym

\centerline{\bf Group-cograded multiplier Hopf (\st)algebras}
\bigskip
\bigskip
\centerline{\it A.T.\ Abd El-hafez
\rm($^{1}$), \it L.\ Delvaux \rm($^{2}$) and \it A.\ Van Daele \rm ($^{3}$)}
\bigskip
\bigskip
\bigskip
{\bf Abstract}
\bigskip
Let $G$ be a group and assume that $(A_p)_{p\in G}$ is a family of
algebras with identity. We have a {\it Hopf\ $G$-coalgebra} (in
the sense of Turaev) if, for each pair $p,q\in G$, there is given
a unital homomorphism $\co_{p,q}:A_{pq}\to A_p \ot A_q$ satisfying
certain properties.
\smallskip
Consider now the direct sum $A$ of these algebras. It is an
algebra, without identity, except when $G$ is a finite group, but
the product is non-degenerate. The maps $\co_{p,q}$ can be used to
define a coproduct $\co$ on $A$ and the conditions imposed on
these maps give that $(A,\co)$ is a multiplier Hopf algebra. It is
$G$-cograded  as explained in this paper.
\smallskip
We study these so-called {\it group-cograded multiplier Hopf
algebras}. They are, as explained above, more general than the
Hopf group-coalgebras as introduced by Turaev. Moreover, our point
of view makes it possible to use results and techniques from the
theory of multiplier Hopf algebras in the study of Hopf
group-coalgebras (and generalizations).
\smallskip
In a separate paper, we treat the quantum double in this context
and we recover, in a simple and natural way (and generalize)
results obtained by Zunino. In this paper, we study integrals, in
general and in the case where the components are
finite-dimensional. Using these ideas, we obtain most of the
results of Virelizier on this subject and consider them in the
framework of multiplier Hopf algebras.

\bigskip
April 2004 ({\it Version 1.0})
\bigskip
\vskip 4 cm
 \hrule
\medskip
($^{1}$) Department of Mathematics, Mansoura University, Mansoura
35516 (Egypt). E-mail: a$_-$t$_-$amer\@yahoo.com
\smallskip
($^{2}$) Department of Mathematics, L.U.C., Universiteitslaan,
B-3590 Diepenbeek (Belgium). E-mail: Lydia.Delvaux\@luc.ac.be
\smallskip
($^{3}$) Department of Mathematics, K.U.\ Leuven, Celestijnenlaan
200B, B-3001 Heverlee (Belgium). E-mail:
Alfons.VanDaele\@wis.kuleuven.ac.be

\newpage

\bf 0. Introduction \rm
\nl
Let $G$ be a (discrete) group. Assume that we have a family of
algebras $(A_p)_{p\in G}$ indexed over $G$. In this paper, we will
 only consider algebras over $\Bbb C$, with or without identity,
but with a non-degenerate product. Consider the {\it direct sum}
$\oplus_{p\in G} A_p$ of these algebras and denote it by $A$.
Elements in $A$ are functions $a$ on $G$ so that $a(p)\in A_p$ for
all $p \in G$ and $a(p)=0$ except for finitely many $p$. Pointwise
operations make $A$ into an algebra. It will not have an identity,
except when all the components have an identity and when the group
$G$ is finite. However, the product will still be non-degenerate.
Therefore, we can consider the multiplier algebra $M(A)$ of $A$
(see [VD1, Definition A.1] for a precise definition). In this
case, $M(A)$ will consist of {\it all} functions $a$ on $G$ with
$a(p)\in M(A_p)$ for all $p\in G$. Recall that $M(A_p)=A_p$ if
$A_p$ has an identity. So, if all components have an identity,
elements in $M(A)$ are functions $a$ on $G$ such that $a(p)\in
A_p$ for all $p$, without further restrictions.
\snl
We will also use $A\ot A$ and its multiplier algebra $M(A\ot A)$.
Elements in $A\ot A$ are functions $a$ on $G\times G$ such that
$a(p,q)\in A_p\ot A_q$ for all $p,q\in G$ and such that $a$ has
finite support. On the other hand, elements in $M(A\ot A)$ are
functions $a$ on $G\times G$ such that $a(p,q)\in M(A_p\ot A_q)$
for all $p,q$, without further restrictions.
\snl
In this paper, we study {\it $G$-cograded multiplier Hopf
(\st)algebras}. This notion will be introduced in Section 1 (cf.\
Definition 1.1). A multiplier Hopf algebra $(A,\co)$ is called
$G$-cograded if $A=\oplus_{p\in G} A_p$ as above and
$$\align \co(A_{pq})(1\ot A_q)&=A_p\ot A_q\\
         (A_p\ot 1)\co(A_{pq})&=A_p\ot A_q
\endalign$$
for all $p,q\in G$. Observe that we consider the algebras $A_p$ as
sitting inside $A$ in the obvious way. When we have a multiplier
Hopf \st algebra, we require the components also to be \st
subalgebras.
\snl
A trivial example is obtained when all algebras $A_p$ are just
$\Bbb C$. In this case $A$ is $K(G)$, the algebra of complex
functions with finite support on $G$ (with pointwise sum and
product). The product in $G$ gives rise to a coproduct $\co_G$ on
$K(G)$ defined by $(\co_G(f))(p,q)=f(pq)$ where $f\in K(G)$ and
$p,q\in G$. Notice that $\co_G$ maps $K(G)$ into the multiplier
algebra $M(K(G)\ot K(G))$ where first $K(G)\ot K(G)$ is identified
with $K(G\times G)$ and then $M(K(G) \ot K(G))$ is identified with
the algebra of all complex functions on $G\times G$. With this
coproduct, $K(G)$ is a multiplier Hopf algebra. It is trivially
$G$-cograded. In fact, it is a $G$-cograded multiplier Hopf \st
algebra if we take the complex conjugate of a function as the
involution.
\snl
Non-trival examples will be given (see e.g.\ Example 1.6), but the
above example is important because of the following result
(Theorem 1.2). It is shown that a multiplier Hopf algebra
$(A,\co)$ is $G$-cograded if and only if there is a non-degenerate
injective homomorphism $\gamma:K(G)\to M(A)$ such that $\gamma$
actually maps into the center of $M(A)$ and moreover
$\co\circ\gamma=(\gamma\ot\gamma)\circ\co_G$. Here, $\co_G$ is the
coproduct on $K(G)$ as defined in the preceding paragraph. We also
need the unique extension of $\gamma\ot \gamma$ from $K(G)\ot
K(G)$ to a homomorphism on $M(K(G)\ot K(G))$ (which exists because
$\gamma$ is assumed to be non-degenerate - see e.g.\ [VD1,
Proposition A.5], and later in this introduction). The proof of
this result is not very hard. Given a $G$-cograded multiplier Hopf
algebra $(A,\co)$, the map $\gamma$ is given by
$(\gamma(f))(p)=f(p)1_p$ where $1_p$ is the identity in the
algebra $M(A_p)$. Conversely, given such a map $\gamma$, the
algebra $A_p$ is defined by $\gamma(f_p)A$ where $f_p$ is the
function with value $1$ on $p$ and $0$ everywhere else. In the
case of a \st algebra, the imbedding $\gamma$ must be a \st
homomorphism.
\snl
As a special case, we get the Hopf group-coalgebras as introduced
by Turaev in [T, Section 11] and studied further by Virelizier in
[V] and others ([Z], [W1] and [W2]). The relation is as follows
(see Theorem 1.5). As in [T], let $G$ be a group and let
$(A_p)_{p\in G}$ be a family of algebras. Now it is assumed that
all these algebras have an identity. Given are also unital
homomorphisms $\co_{p,q}:A_{pq}\to A_p\ot A_q$ for all $p,q\in G$.
It is assumed that the family is coassociative in the sense that
$$(\co_{p,q}\ot\iota)\co_{pq,r}=(\iota\ot\co_{q,r})\co_{p,qr}$$
on $A_{pqr}$ for all triples $p,q,r\in G$. We use $\iota$ to
denote the identity map on all $A_p$. If there is also a {\it
counit} and an {\it antipode} (see Theorem 1.5 for a precise
definiton), this system is called a Hopf $G$-coalgebra (see [T,
Section 11]).
\snl
Now assume that we have such a Hopf $G$-coalgebra. Let $A$ be the
direct sum of the algebras $(A_p)_{p\in G}$ as before. Define
$\co:A\to M(A\ot A)$ by
$$(\co(a))(p,q)=\co_{p,q}(a(pq))$$
when $a\in A$ for all $p,q\in G$. Recall that $a(pq)\in A_{pq}$
and that $\co_{p,q}:A_{pq}\to A_p\ot A_q$ so that the right hand
side of the above equation indeed belongs to $A_p\ot A_q$ and
therefore, $\co(a)$ is really defined in $M(A\ot A)$ by the above
formula. It is easy to see that $\co$ is a non-degenerate
homomorphism from $A$ in $M(A\ot A)$. It is coassociative because
the family $(\co_{p,q})_{p,q}$ is assumed to be coassociative.
\snl
It will be shown in Theorem 1.5 that $(A,\co)$, as defined above,
is actually a multiplier Hopf algebra (in the sense of [VD1]).
This is so because of the existence of the counit and the antipode
above and it is proven in a more or less standard way. It is also
shown that the original direct sum decomposition makes it into a
$G$-cograded multiplier Hopf algebra. In general however, not
every $G$-cograded multiplier Hopf algebra $A$ comes from a Hopf
$G$-coalgebra. For this to happen, it is necessary and sufficient
that $A$ is regular and that all the components are algebras with
identity. All this is proven in Theorem 1.5.
\snl
We see that, roughly speaking, the Hopf $G$-coalgebras, introduced
by Turaev and studied further by others, are in fact nothing else
but multiplier Hopf algebras $(A,\co)$ with a $G$-cograding, or
equvalently, with a distinguished embedding of $K(G)$ into the
center of $M(A)$.
\snl
This point of view is important. It allows the results and
techniques from the theory of multiplier Hopf algebras to be used
in the study of Hopf group-coalgebras (and its generalizations).
\snl
In a separate paper ([D-VD]), we follow these ideas to study the
quantum double for Hopf group-coalgebras. We generalize the
setting, considered in the paper by Zunino [Z], and construct a
family of associated multiplier Hopf algebras. One of them is the
usual Drinfel'd double, the other one is the double constructed in
[Z].
\snl
In this paper however, we mainly study integrals on group-cograded
multiplier Hopf algebras and we obtain our results by applying
known result about multiplier Hopf algebras in general. We obtain
most of the results by Virelizier on this subject ([V]), in a
simpler, natural way (and in greater generality).
\nl
Let us briefly summarize the {\it content} of the paper.
\snl
In {\it Section 1} we introduce the notion of a group-cograded
multiplier Hopf (\st)algebra (as explained earlier in this
introduction). We focus on the relation with the Hopf
group-coalgebras, introduced by Turaev. In {\it Section 2}, we
apply the known results about integrals on multiplier Hopf
algebras to our group-cograded multiplier Hopf algebras. This is a
small section because the results are quite obvious and easy to
obtain. Finally, in the last section, {\it Section 3}, we fully
make use of known techniques in the theory of multiplier Hopf
algebras, more precisely about multiplier Hopf algebras of
discrete type, to construct integrals on a group-cograded
multiplier Hopf algebra when the components are all
finite-dimensional. We recover results obtained by Virelizier [V].
We also give an explicit formula for the integrals.
\nl
The main material, needed for reading this paper, is found in the
basic references on (regular) multiplier Hopf (\st)algebras
([VD1], [VD-Z2]) and multiplier Hopf algebras with integrals
[VD4]. We will freely use notions and results of these  papers.
However, if convenient, we will recall the main concepts so as to
make the paper self-contained to a certain extend. We will also
use the Sweedler notation in the case of multiplier Hopf algebras.
It is known that this is justified and often, it makes formulas
and arguments more easy to understand. For the general results on
Hopf algebras, we refer to the standard works of Abe [A] and
Sweedler [S].
\nl
Finally, let us recall some of the {\it conventions} and {\it
notations} that will be used. We use $1$ for the identity in
various algebras while we use $e$ for the identity element in a
group. The symbol $\iota$ will be reserved for the identity map,
usually here on algebras and subalgebras. We will not give these
symbols an index to indicate e.g.\ in what algebra we have the
identity, except when it is really necessary to avoid confusion.
Usually however, things should be clear from the context.
Similarly, we will always use $\co$ to denote a comultiplication.
\snl
Because this plays an important role, we will say something more
about non-degenerate homomorphisms and their extensions. Take two
algebras $A$ and $B$, both with or without identity, but wit a
non-degenerate product. Consider their multiplier algebras $M(A)$
and $M(B)$ respectively. Let $\alpha:A \to M(B)$ be a homomorphism
(i.e.\ an algebra map). It is called non-degenerate if
$\alpha(A)B=B\alpha(A)=B$. Then, it is possible to extend $\alpha$
to a unital homomorphism from $M(A)$ to $M(B)$. This extension is
uniquely defined and therefore, it makes sense to denote also the
extension with the same symbol $\alpha$. This can e.g.\ be applied
to the comultiplication $\co$ on a multiplier Hopf algebra $A$.
The extension is a unital homomorphism from $M(A)$ to $M(A\ot A)$.
It can also be applied to the counit. Then the extension becomes a
homomorphism from $M(A)$ to $\Bbb C$. Finally, also the antipode
$S$ of a (possibly non-regular) multiplier Hopf algebra can be
extended to a map from $M(A)$ to itself. In this case one has to
apply the general result with $B=A$ as a vector space, but endowed
with the opposite multiplication so as to make $S:A\to M(B)$ a
homomorphism. That the antipode $S$ is non-degenerate follows
easily from the formulas, characterizing the antipode, given by
$$\align  m((S\ot\iota)(\co(a)(1\ot b)))&=\varepsilon(a)b \\
m((\iota\ot S)((b\ot 1)\co(a))) &=\varepsilon(a)b
\endalign$$
for $a,b \in A$. In these formulas, we have used $m$ to denote the
multiplication in $A$, as a map from $A\ot A$ to $A$, and extended
to a map on $M(A)\ot A$ and on $A\ot M(A)$. For properties of
non-degenerate homomorphisms and extensions, we refer to the
appendix of [VD1].

\nl\nl
\bf Acknowlegdements \rm
\smallskip
We like to thank S.\ Wang for fruitful discussions on this
subject. The first author wants to thank Prof.\ A.S.\ Hegazi for
his continuous encouragement and M.\ Mansour for helpful
discussions.

\nl\nl

\bf 1. Group-cograded multiplier Hopf (\st)algebras \rm
\nl
Let $G$ be any (discrete) group. We now introduce the notion of a
$G$-cograded multiplier Hopf (\st)algebra. It is the main object
in this paper.

\inspr{1.1} Definition \rm
Let $(A,\co)$ be a multiplier Hopf algebra. Assume that there is a
family of (non-trivial) subalgebras $(A_p)_{p\in G}$ of $A$ so
that
\item{} i) $A=\sum_{p\in G} \oplus_{p\in G} A_p$ with $A_pA_q=0$ whenever
$p,q\in G$ and $p\neq q$,
\item{} ii) $\co(A_{pq})(1\ot A_q)=A_p\ot A_q$ and $(A_p\ot
1)\co(A_{pq})=A_p\ot A_q$ for all $p,q\in G$.
\snl
Then we call $(A,\co)$ a {\it $G$-cograded} multiplier Hopf
algebra.
\einspr

If $(A,\co)$ is a multiplier Hopf \st algebra, we will assume that
all these subalgebras are \- \st subalgebras and call it a
$G$-cograded multiplier Hopf \st algebra. Examples of $G$-cograded
multiplier Hopf algebras will be given later in this section (see
Example 1.6).
\snl
In the following theorem we characterize a $G$-cograded multiplier
Hopf algebra by using the multiplier Hopf algebra $(K(G),\co)$ as
defined in the introduction.

\inspr{1.2} Theorem \rm
A multiplier Hopf algebra $(A,\co)$ is $G$-cograded if and only if
there is a non-degenerate algebra embedding $\gamma:K(G)\to M(A)$
satisfying
\item{} i) $\gamma(K(G))$ belongs to the centre $Z(M(A))$ of
$M(A)$,
\item{} ii) $\gamma$ is also a coalgebra map, that is,
$\co(\gamma(f))=(\gamma\ot \gamma)(\co(f))$ for all $f\in K(G)$.
\snl
In the case of a $G$-cograded multiplier Hopf \st algebra, we have
that $\gamma$ is a \st map.
\snl \bf Proof: \rm
First assume that $A$ is a $G$-cograded multiplier Hopf algebra.
We write $A=\oplus_{p\in G} A_p$ as in Definition 1.1. It is not
hard to see that $M(A)=\prod_{p\in G}M(A_p)$. Define a map
$\gamma:K(G)\to M(A)$ by $(\gamma(f))(p)=f(p)1_p$ for any $p\in G$
where $1_p$ is the identity in $M(A_p)$. It is easy to see that
$\gamma$ is a homomorphism and that it maps into the centre of
$M(A)$.
\snl
To prove that it is non-degenerate, take any $a\in A$. Because $A$
is the direct sum of the algebras $A_p$, we have a finite subset
$F$ of $G$ such that $a(p)=0$ for $p\notin F$. If $f\in K(G)$ is
such that $f(p)=1$ when $p\in F$, we have that $\gamma(f)a=a$.
Therefore, $\gamma$ is a non-degenerate homomorphism from $K(G)$
to $M(A)$.
\snl
For any $f\in K(G)$ and $p,q\in G$ we get
$$((\gamma\ot\gamma)(\co(f)))(p,q)=(\co(f))(p,q)(1_p\ot
   1_q)=f(pq)\co(1_{pq})$$
and we see that $\co(\gamma(f))=(\gamma\ot \gamma)(\co(f))$. This
completes the proof of one direction.
\snl
Conversely, suppose that we have a multiplier Hopf algebra
$(A,\co)$ with a non-degener\-ate homomorphism $\gamma:K(G)\to
M(A)$ satisfying the conditions of the theorem. Define the
subalgebra $A_p$ by $\gamma(f_p)A$ where $p\in G$ and $f_p$ is the
function on $G$ satisfying $f_p(p)=1$ and $f_p(q)=0$ when $q\neq
p$. Because $\gamma$ maps into the centre, $A_p$ will be a
subalgebra of $A$. It is non-trivial for each $p$ because $\gamma$
is assumed to be injective. The product in each subalgebra is
still non-degenerate. And because $\gamma$ is assumed to be
non-degenerate, we have $A=\gamma(K(G))A$ and it follows that
$A=\oplus_{p\in G} A_p$. Clearly also $A_pA_q={0}$ if $p\neq q$.
The identity $1_p$ in $M(A_p)$ is $\gamma(f_p)$.
\snl
From the fact that $\gamma$ respects the comultiplications in the
sense that $\co(\gamma(f))=(\gamma\ot \gamma)(\co(f))$ for all
$f\in K(G)$, we get
$$\align \co(A_{pq})(1\ot A_q)&=\co(\gamma(f_{pq})A)(1\ot
\gamma(f_q)A))\\
&=(\gamma\ot\gamma)(\co(f_{pq})(1\ot f_q))(\co(A)(1\ot A))\\
&=((\gamma\ot\gamma)(f_p\ot f_q))(A\ot A)=A_p\ot A_q \endalign$$
for all $p,q\in G$. Similarly, the other equality is proven. This
completes the proof of the other implication in the theorem.
\einspr

We will use this charaterization to look at the counit and the
antipode of a group-cograded multiplier Hopf algebra. First we
prove the following lemma. We use, as explained already in the
introduction, extensions to the multiplier algebra of
non-degenerate homomorphisms and use the same symbol to denote
these extension.

\inspr{1.3} Lemma \rm
Let $A$ and $B$ be multiplier Hopf algebras. Assume that
$\alpha:A\to M(B)$ is a non-degenerate homomorphism that respects
the comultiplication in the sense that $\co(\alpha(a))=(\alpha\ot
\alpha)\co(a)$ for all $a\in A$. Then $\alpha$ preserves the unit,
the counit and the antipode, that is,
$$\alpha(1)=1 \qquad\qquad \varepsilon(\alpha(a))=\varepsilon(a)
\qquad\qquad S(\alpha(a))=\alpha(S(a))$$ for all $a\in A$.

\snl\bf Proof: \rm
As $\alpha$ is a non-degenerate homomorphism from $A\to M(B)$, the
extension is a unital homomorphism from $M(A)$ to $M(B)$ as we
have said already in the introduction.
\snl
To show that $\alpha$ respects the counits, take any $a,a'\in A$
and $b\in B$. We can apply the counit in $B$ on the first leg of
the equation
$$\co(\alpha(a))(1\ot \alpha(a')b)=((\alpha\ot\alpha)(\co(a)(1\ot a')))(1\ot b)$$
to get
$\alpha(aa')b=\sum_{(a)}\varepsilon(\alpha(a_{(1)}))\alpha(a_{(2)}a')b$.
On the other hand, we also have
$aa'=\sum_{(a)}\varepsilon(a_{(1)})a_{(2)}a'$. Then we use the
fact that $\co(A)(1\ot A)=A\ot A$ and we can conclude that
$\varepsilon(c)\alpha(c')b=\varepsilon(\alpha(c))\alpha(c')b$ for
all $c,c'\in A$ and $b\in B$. Now because $\alpha(A)B=B$ we get
that $\varepsilon(c)=\varepsilon(\alpha(c))$ for all $c\in A$.
\snl
Next we show that $\alpha$ converts one antipode to the other.
Again take any $a,a'\in A$ and $b\in B$ and consider the equation
$$\co(\alpha(a))(1\ot \alpha(a')b)=((\alpha\ot\alpha)(\co(a)(1\ot a')))(1\ot
  b).$$
Now, we apply the map $m\circ(S\ot \iota)$ on both sides of this
equation. Notice that we have elements in $M(B)\ot B$ and that it
is also possible to extend the antipode to $M(B)$ (as we mentioned
in the introduction). We get
$$\varepsilon(\alpha(a))\alpha(a')b
          =\sum_{(a)}S(\alpha(a_{(1)}))\alpha(a_{(2)}a')b.$$
Because $\varepsilon(\alpha(a)=\varepsilon(a)$, we obtain
$$\sum_{(a)}\alpha(S(a_{(1)}))\alpha(a_{(2)}a')b=
   \sum_{(a)}S(\alpha(a_{(1)}))\alpha(a_{(2)}a')b.$$
Again, because $\co(A)(1\ot A)=A\ot A$ we conclude that
$\alpha(S(c))\alpha(c')b= S(\alpha(c))\alpha(c')b$ for all
$c,c'\in A$ and $b\in B$. Using once more that $\alpha$ is
non-degenerate, we get $\alpha(S(c))=S(\alpha(c))$ for all $c\in
A$. This completes the proof.
\einspr

The reader may notice that the above arguments are very similar to
the onces used to show that the counit and the antipode in
multiplier Hopf algebras are unique. This is not surprising. It
would be a consequence of the result in the lemma that the counit
and the antipode is unique. Simply take $A=B$ with the identity
map for $\alpha$.
\snl
Remark that we do not need the image of $\alpha$ to belong to the
center for the above lemma.
\snl
It is now easy to obtain the following result.

\inspr{1.4} Proposition \rm
Let $A$ be a $G$-cograded multiplier Hopf algebra with
decomposition $A=\oplus_{p\in G} A_p$ as in Definition 1.1. Then
$\varepsilon(A_p)=0$ when $p\neq e$ and $S(A_p)\subseteq
M(A_{p^{-1}})$ for all $p$. If the multiplier Hopf algebra is
regular, we have $S(A_p)=A_{p^{-1}}$.

\snl\bf Proof: \rm
By Theorem 1.2 there is a non-degenerate central imbedding
$\gamma:K(G)\to M(A)$ which respects the comultiplications. So
Lemma 1.3 applies and the result follows from the fact that the
counit and the antipode on $K(G)$ are given by the formulas
$\varepsilon(f)=f(e)$ and $(S(f))(p)=f(p^{-1})$ for $f\in K(G)$
and $p\in G$.
\einspr

In the following theorem, we compare this notion with the Hopf
group-coalgebras as introduced by Turaev in [T, Section 11]. It
shows that, roughly speaking, Hopf group-coalgebras are multiplier
Hopf algebras with some extra structure. This observation throws a
new light on the notion introduced by Turaev. We will see some of
the consequences in the rest of the paper (see Section 2 and 3).
We also refer to our work on the quantum double in this context
[D-VD].

\inspr{1.5} Theorem \rm
Let $G$ be a group. Assume that for all $p\in G$, we have given an
algebra $A_p$, with or without identity, but with a non-degenerate
product. Moreover, we require the following:
\snl
i) For all $p,q\in G$ we have a homomorphism $\co_{p,q}:A_{pq}\to
M(A_p\ot A_q)$ satisfying $\co_{p,q}(A_{pq})(1\ot A_q)\subseteq
A_p\ot A_q$ and $(A_p\ot 1)\co_{p,q}(A_{pq})\subseteq A_p\ot A_q$,
as well as $(1\ot A_q)\co_{p,q}(A_{pq})\subseteq A_p\ot A_q$ and
$\co_{p,q}(A_{pq})(A_p\ot 1)\subseteq A_p\ot A_q$. Furthermore, it
is assumed that these maps form a 'coassociative' family in the
sense that for all $p,q,r\in G$ we have
$$(c\ot 1\ot 1)((\co_{p,q}\ot \iota)(\co_{pq,r}(a)(1\ot b))) =
((\iota\ot\co_{q,r})((c\ot 1)\co_{p,qr}(a)))(1\ot 1\ot b)$$
whenever $a\in A_{pqr}$, $b\in A_r$ and $c\in A_p$.
\snl
ii) There is a homomorphism $\varepsilon_e:A_e\to \Bbb C$ so that
for all $p\in G$ we get
$$\align(\iota\ot \varepsilon_e)((a\ot 1)\co_{p,e}(b)) &= ab\\
  (\varepsilon_e\ot \iota)(\co_{e,p}(a)(1\ot b)) &=ab
\endalign$$
whenever $a,b\in A_p$.
\snl
iii) For all $p\in G$, there is a anti-isomorphism $S:A_p\to
A_{p^{-1}}$ so that
$$\align m(S_{p^{-1}}\ot \iota)(\co_{p^{-1},p}(a)(1\ot b))
&=\varepsilon_e(a)b\\
m(\iota \ot S_{p^{-1}})((b\ot 1)\co_{p,p^{-1}}(a))
&=\varepsilon_e(a)b
\endalign
$$
for all $p\in A_e$ and $b\in A_p$
\snl
Under these assumptions, the algebra $A$, defined as $\oplus_{p\in
G} A_p$, can be given a comultiplication in a natural way, making
it into a regular multiplier Hopf algebra.
\snl
Conversely, assume that we have a regular $G$-cograded multiplier
Hopf algebra $A$ with decomposition $A=\oplus_{p\in G} A_p$ as in
Definition 1.1. Then, the family of subalgebras $(A_p)_{p\in G}$
can be endowed (in a natural way) with the above objects and they
will satisfy the above conditions i), ii) and iii).

\snl \bf Proof: \rm
First, assume that we have such a family of algebras, provided
with these data satisfying the condition i), ii) and iii) above.
Define the algebra $A=\oplus_{p\in G} A_p$ as before. Recall that
this is an algebra with non-degenerate product, that hence we can
consider the multiplier algebra $M(A)$ and that this is given by
the product $\prod_{p\in G}M(A_p)$. Also $A\ot A=\oplus_{p,q}
(A_p\ot A_q)$ while $M(A\ot A)=\prod_{p,q}M(A_p\ot A_q)$. We can
define a map $\co:A\to M(A\ot A)$ by
$$(\co(a))(p,q)=\co_{p,q}(a(pq))$$
whenever $p,q\in G$. Because $a(pq)\in A_{pq}$ and $\co_{p,q}$
maps $A_{p,q}$ into $M(A_p\ot A_q)$, this map is well-defined.
\snl
Because of all the conditions on these maps $\co_{p,q}$, it
follows easily that $\co(A)(1\ot A)$ and $(1\ot A)\co(A)$, as well
as $\co(A)(A\ot 1)$ and $(1\ot A)\co(A)$, all actually are subsets
of $A\ot A$. Because the maps $\co_{p,q}$ are assumed to be
homomorphisms, the same is true for $\co$. Finally,
coassociativity of $\co$ follows by the 'coassociativity' of the
family $(\co_{p,q})_{p,q}$.
\snl
We will now use Proposition 2.9 from [VD4] to show that the pair
$(A,\co)$ is indeed a (regular) multiplier Hopf algebra. To use
this proposition, we have to define the counit and the antipode on
$A$. We first define $\varepsilon$ on $A$ by
$\varepsilon(a)=\varepsilon_e(a(e))$. It is straigthforward to
verify that it is a homomorphism and satisfies the requirements of
a counit. The antipode $S$ is defined on $A$ by
$(S(a))(p^{-1})=S_p(a(p))$ for all $p\in G$. Again it is
straightforward to show that this is a anti-homomorphism from $A$
to itself and that it satisfies the requirements of an antipode.
Hence, by Proposition 2.9 of [VD4], we do have a regular
multiplier Hopf algebra. It is clearly $G$-cograded in the sense
of Definition 1.1.
\snl
Conversely, now assume that we have a $G$-cograded multiplier Hopf
algebra $A$ with $A=\oplus_{p\in G} A_p$ as in the definition. We
now define $\co_{p,q}:A_{pq}\to M(A_p\ot A_q)$ by the restricting
$\co$ to $A_{pq}$. These maps will satisfy the requirements.
Moreover, because of Theorem 1.2 and Lemma 1.3, we must have that
the counit $\varepsilon$ and antipode $S$ satifsy
$\varepsilon(a)=0$ if $a\in A_p$ and $p\neq e$, as well as
$S(a)\in A_{p^{-1}}$ when $a\in A_p$. Then we can define
$\varepsilon_e$ on $A_e$ simply by restricting $\varepsilon$ to
$A_e$ and we can define $S_p$ on $A_p$  also by taking the
restriction of $S$. Again it is all straightforward to verify that
the assumptions in i) and ii) are satisfied.
\snl
This proves the result.
\einspr

As a special case of the above theorem, consider a Hopf
group-coalgebra as defined by Turaev in [T]. This is exactly the
situation of Theorem 1.5 with all the components unital algebras.
\snl
Before we pass to examples, we would like to make some comment on
the notion of regularity and why it is needed for the above
result. Recall first that a multiplier Hopf algebra $(A,\co)$ is
called regular if also $(A,\co^{\text{op}})$ is a multiplier Hopf
algebra (see Definition 2.3 in [VD1]). This is exactly the case
when the antipode $S$ maps $A$ into itself (and not only in $M(A)$
as it does in general) and when it is bijective. This is used in
the previous proof because we refer to Proposition 2.9 in [VD4] to
show that we actually get a multiplier Hopf algebra. Probably, it
is possible to formulate some similar result that holds also when
we do not have regularity.
\snl
On the other hand, for the Hopf group-coalgebras of Turaev, we
know that the antipodes are assumed to be bijections. And also, in
the case of a Hopf \st algebra, regularity is automatic (see again
[VD1]).
\nl
We finish this section by looking at some examples.
\snl
We have already considered the basic but trivial example where the
multiplier Hopf algebra is $K(G)$, the algebra of functions with
finite support on a group $G$, and all the components just the
trivial algebra $\Bbb C$. It is not hard to use this to construct
other, still rather trivial examples. Indeed, take any multplier
Hopf algebra $A_0$ and any group $G$. Now let $A=K(G)\ot A_0$ with
the tensor product algebra structure and the tensor product
coalgebra structure. So, $(f\ot a)(g\ot b)=fg\ot ab$ when $f,g\in
K(G)$ and $a,b\in A_0$. Also $\co(f\ot
a)=(\iota\ot\sigma\ot\iota)(\co(f)\ot \co(a))$ where $f\in K(G)$
and $a\in A_0$ and where $\sigma$ is the flip from $K(G)\ot A_0$
to $A_0\ot K(G)$. In this case, the imbedding $\gamma:K(G)\to
M(A)$ is of course given by $\gamma(f)=f\ot 1$ where $1$ is the
identity in $M(A_0)$. All the components are $A_0$ and all the
maps $\co_{p,q}$ are simply the comultiplication on $A_0$.
\snl
A non-trivial example is obtained from the above by deforming the
comultiplication with the help of an action of the group on $A$.
This gives the following non-trivial example. It is known as a
multiplier Hopf algebra (see e.g.\ Example 3.3 in [D]), but we
present it here as an example of a group-cograded multiplier Hopf
algebra.

\inspr{1.6} Example \rm
Let $(A_0,\co_0)$ be a multiplier Hopf algebra. Let $G$ be a group
and assume that $G$ acts on the multiplier Hopf algebra $A_0$ in
the following sense. For each $p\in G$, we have an automorphism
$\alpha_p$ of the algebra $A_0$ such that also
$\co_0(\alpha_p(a))=(\alpha_p\ot\alpha_p)\co_0(a)$ for all $a\in
A_0$ and all $p\in G$. We also assume that $\alpha_e=\iota$ and
$\alpha_p(\alpha_q(a))=\alpha_{pq}(a)$ for all $a\in A_0$ and
$p,q\in G$.
\snl
Now, we consider again the algebra $A$, defined as $K(G)\ot A_0$
with the tensor product algebra structure. We will consider
elements in $A$ as functions on $G$ with finite support and values
in $A_0$. This algebra is made into a multiplier Hopf algebra by
defining the coproduct $\co$ by
$$(\co(a))(p,q)=(\alpha_q\ot\iota)(\co_0(a(pq)))$$
where $a\in A$ and $p,q\in G$. It is straightforward to show that
this is a multiplier Hopf algebra. The counit $\varepsilon$ on $A$
is given by $\varepsilon(a)=\varepsilon_0(a(e))$ where
$\varepsilon_0$ is the counit on $A_0$. The antipode $S$ on $A$ is
given by $S(a)(p)=S_0(\alpha_p(a(p^{-1})))$. If the action is
trivial, we get the tensor product coalgebra structure as
discussed above.
\snl
Now, we have the imbedding $\gamma:K(G)\to M(A)$ given by
$\gamma(f)=f\ot 1$ where again $1$ stands for the identity in
$M(A_0)$. Of course, this is a central imbedding and it is trivial
to verify that it is compatible with the comultiplications. So, we
have a G-cograded multiplier Hopf algebra. Again, all the
components are the same, namely $A_0$. However, not all the maps
$\co_{p,q}$ are equal. We have $\co_{p,q}=(\alpha_q\ot\iota)\co_0$
for all $p,q\in G$.
\einspr

It is not hard to find more concrete examples using the above
construction. Take e.g. another group $H$ and assume that $G$ acts
on $H$ by means of automorphisms $\rho_p$, with $p\in G$. Take
$A_0=K(H)$ with its comultiplication coming from the product in
the group $H$. Then an action $\alpha$ of $G$ on $A_0$ is given by
the formula $(\alpha_p(f))(h)=f(\rho_{p^{-1}}(h))$ for $f\in K(H)$
and $h\in H$ and $p\in G$. It is not hard to see that we get the
algebra of complex functions with finite support on the
semi-direct product of $G$ with $H$.
\snl
One can still make it more concrete if we let $H=G$ and take the
adjoint action $\rho$ of $G$ on itself given by
$\rho_p(q)=pqp^{-1}$ when $p,q\in G$. Then $A$ is $K(G\times G)$
while the coproduct is given by
$$(\co(f))(p,h,q,k)=f(pq,q^{-1}hqk)$$
whenever $p,q,k,h\in G$. The embedding $\gamma$ is now given by
$(\gamma(f))(p,h)=f(p)$.
\snl
It is also possible to give examples of related situations.
Consider e.g. the discrete quantum group $su_q(2)$, the Pontryagin
dual of the compact quantum group $SU_q(2)$. As an algebra, it is
the direct sum $A$ of the matrix algebras $M_n(\Bbb C)$ where
$n=1,2,3,\ldots$. The comultiplication is rather complicated on
this level. Inside the multiplier algebra $M(A)$, we have the
elements $k,e,f$ satisfying the commutation rules $ke=\lambda ek$,
$kf=\lambda^{-1}fk$ and $ef-fe=k^{2}-k^{-2}$. Here $\lambda$ is
the deformation parameter and $0<\lambda<1$. It is possible to
construct an imbedding $\gamma$ of $K(\Bbb Z)$ into the multiplier
algebra $M(A)$ of this discrete quantum group. It sends the delta
function $f_p$, defined as $1$ in $p$ and $0$ anywhere else, to
the spectral projection of $k$ associated with the eigenvalue
$\lambda^{\frac12 p}$. This imbedding respects the
comultiplication. However, it is not central as $k$ is not a
central element in $A$.
\snl
This example shows that more general objects than the
group-cograded multiplier Hopf algebras are also of interest. The
requirement that $K(G)$ is imbedded in the center of $M(A)$ may be
somewhat restrictive. In [L-VD], related objects are studied
without this restriction so that the above example (with the
discrete quantum group $su_q(2)$) fits into the theory.
\nl\nl

\bf 2. Integrals on group-cograded multiplier Hopf (\st)algebras
\rm
\nl
In this section, we consider a regular multiplier Hopf algebra
$(A,\co)$. Regular multiplier Hopf algebras with integrals are
studied in [VD4]. We recall some of the results.
\snl
First observe that for any linear functional $\varphi$ on $A$ and
any element $a\in A$, we can define a multiplier
$(\iota\ot\varphi)\co(a)$ in $M(A)$ by the formulas
$$\align ((\iota\ot\varphi)\co(a))b &=
(\iota\ot\varphi)(\co(a)(b\ot 1)) \\
b((\iota\ot\varphi)\co(a)) &= (\iota\ot\varphi)((b\ot 1)\co(a))
\endalign$$
where $b\in A$. Because we assume our multiplier Hopf algebra to
be regular, we have not only $(b\ot 1)\co(a)\subseteq A\ot A$ but
also $\co(a)(b\ot 1)\subseteq A\ot A$ and so the above formulas
make sense.

\inspr{2.1} Definition \rm
A linear functional $\varphi$ on $A$ is called {\it left
invariant} if $(\iota\ot\varphi)\co(a)=\varphi(a)1$ in $M(A)$ for
all $a\in A$. A non-zero left invariant functional is called a
{\it left integral} on $A$. Similarly, a non-zero linear function
$\psi$ satisfying $(\psi\ot\iota)\co(a)=\psi(a)1$ for all $a\in A$
is called a {\it right integral} on $A$.
\einspr

There are various data (with many relations among them) associated
with left and right integrals. We recall the following definition
and results from [VD4].

\inspr{2.2} Theorem \rm
Let $(A,\co)$ be a regular multiplier Hopf algebra with a left
integral $\varphi$. Any other left integral is a scalar multiple
of $\varphi$. There exists also a right integral, unique up to a
scalar, given by $\varphi\circ S$. There is a scalar $\nu$ (the
{\it scaling constant}) in $\Bbb C$ defined by $\varphi\circ
S^2=\nu\varphi$. There is a grouplike multiplier $\delta$ in
$M(A)$ (the {\it modular element}) such that
$(\varphi\ot\iota)\co(a)=\varphi(a)\delta$ for all $a\in A$. The
integral is faithful in the sense that for any $a\in A$ we have
$a=0$ if either $\varphi(ab)=0$ for all $b\in A$ or
$\varphi(ba)=0$ for all $b$. Also the right integral is faithful.
Finally, there are automorphisms $\sigma$ and $\sigma'$ (the {\it
modular automorphims}) satisfying
$\varphi(ab)=\varphi(b\sigma(a))$ and $\psi(ab)=\psi(b\sigma'(a))$
for all $a,b\in A$.

\inspr{2.3} Definition \rm
As before, let $(A,\co)$ be a regular multiplier Hopf algebra with
a left integral $\varphi$. Denote by $\hat A$ the space of linear
functionals on $A$ of the form $x\mapsto \varphi(xa)$ for some
$a\in A$. Then $\hat A$ is made into a regular multiplier Hopf
algebra with the product and coproduct in $\hat A$ dual to the
coproduct and product resp.\ in $A$. This is called the {\it dual}
of $A$.
\einspr

It is again a multiplier Hopf algebra with integrals. Moreover,
the dual of $\hat A$ is canonically isomorphic with the original
multiplier Hopf algebra. This is some form of the Pontryagin
duality for multiplier Hopf algebras with integrals. It is a
generalization of the duality theorem for finite-dimensional Hopf
algebras.
\snl
We have two important special cases. First, we have the multiplier
Hopf algebras of {\it compact type}. These are nothing else but
Hopf algebras with integrals. The other ones are the multiplier
Hopf algebras of {\it discrete type}. A multiplier Hopf algebra is
called of discrete type if it has a co-integral. A {\it left
co-integral} is a non-zero element $h\in A$ satisfying
$ah=\varepsilon(a)h$ for all $a\in A$.
\snl
In [VD-Z1, Theorem 2.10], it is proven that there always exists an
integral on a multiplier Hopf algebra of discrete type. Moreover,
the dual $\hat A$ of a multiplier Hopf algebra of discrete type is
of compact type, that is a Hopf algebra (with integrals); see
again Proposition 2.1 in [VD-Z1]. We will use some of these
results in the next section.
\nl
We now will apply the above results about integrals to regular
group-cograded multiplier Hopf algebras as they were introduced in
the previous section.
\snl
Let $G$ be a group and $(A,\co)$  a regular $G$-cograded
multiplier Hopf algebra with decomposition $A=\oplus_{p\in G}
A_p$. The (full) linear dual $A'$ of $A$ is given by $\prod_{p\in
G}(A_p)'$. Here elements are given by functions $f$ on $G$ such
that $f(p)\in A_p'$ for all $p\in G$ (with no restrictions on the
support). We will write $f_p$ in stead of $f(p)$. The multiplier
algebra $M(A)$ is given by $\prod_{p\in G}M(A_p)$. Again, when
$m\in M(A)$, we will consider the associated elements $m_p$ in
$M(A_p)$. Let $\co_{p,q}$ be the associated map from $A_{pq}$ to
$M(A_p\ot A_q)$, as given in the proof of Proposition 1.5.
\snl
The following results are easy to obtain.

\inspr{2.4} Proposition \rm
Let $A$ be a regular $G$-cograded multiplier Hopf algebra as
before. Let $\varphi$ be a left integral on $A$. Then all the
components $\varphi_p$ are faithful. We have
$(\iota\ot\varphi_q)(\co_{p,q}(a))=\varphi_{pq}(a)1_p$ for all
$a\in A_{pq}$ and all $p,q\in G$.

\inspr{2.5} Proposition \rm
Let $A$ and $\varphi$ be as in the previous proposition. For the
components $(\delta_p)$ of the modular element $\delta$, we have
$(\varphi_p\ot \iota)\co_{p,q}(a) = \varphi_{pq}(a)\delta_q$ when
$a\in A_{pq}$ and $p,q\in G$. Also $\co_{p,q}(\delta_{pq})(b\ot
c)=\delta_p b\ot \delta_q c$ for all  $a\in A_{pq}$, $b\in A_p$,
$c\in A_q$ and $p,q\in G$.
\einspr

There is not much one can say about the scaling constant, except
that it is the same for all components. Indeed, as the square of
the antipode leaves each component globally invariant, we have
trivially $\varphi_p(S^2(a))=\nu \varphi_p(a)$ whenever $a\in A_p$
and $p\in G$.
\snl
Also, when we look at the modular automorphisms, we see that they
again must leave the components globally invariant. If we denote
by $\sigma_p$ the restriction of $\sigma$ to $A_p$ we have
$\varphi_p(ab)=\varphi_p(b\sigma_p(a))$ for all $a,b\in A_p$ and
all $p\in G$.
\snl
In fact, all the relations among the different data associated
with a regular multiplier Hopf algebra $A$ with integrals (a
so-called algebraic quantum group), like $\delta, \sigma, \sigma'$
and their dual objects related with the dual $\hat A$, will give
rise to similar relations among the components.
\nl
Let us finish by looking briefly at the example in the previous
section (Example 1.6). Suppose that the multiplier Hopf algebra
$A_0$ has a left integral $\varphi_0$. By the uniqueness of
integrals, we must have a homomorphism $\mu:G\to \Bbb C$, given by
$\varphi(\alpha_p(a))=\mu(p)\varphi(a)$ for all $a\in A$ and $p\in
G$. The new multiplier Hopf algebra $A$ will also have a left
integral $\varphi$ and it will be given by $\varphi(a)=\sum_{p\in
G}\varphi_0(a(p))$ where $a\in A$ and where we consider $a$ as a
function on $G$ with finite support and with values in $A_0$ as in
Example 1.6. For the modular element $\delta$ in $M(A)$ we get
components $\delta_p=\mu(p)\delta_0$ for all $p$.

\nl\nl

\bf 3. Group-cograded multiplier Hopf (\st) algebras with
finite-dimensional components \rm
\nl
In this section, we consider regular group-cograded multiplier
Hopf algebras with finite-dimensional components. So, as before,
let $G$ be a group and let $(A,\co)$ be a $G$-cograded multiplier
Hopf algebra. Let $A=\oplus_{p\in G} A_p$ as in the definition and
assume now that every subalgebra $A_p$ has a unit and is
finite-dimensional. As we have seen in Section 1, we then have a
Hopf $G$-coalgebra. It is called 'finite-dimensional' by Turaev in
[T, Section 11].
\snl
We will now use results from the general theory of multiplier Hopf
algebras to obtain properties in this special case. We will show
that in this case, the multiplier Hopf algebra $A$ is of discrete
type in the sense of [VD-Z1]. It then follows that there are
integrals on $A$. We recover the results of Virelizier obtained in
[V]. Moreover, we give an explicit formula for the left integral
$\varphi$ in terms of a basis and a dual basis in $A_p$ and its
dual $A_p'$ for each $p\in G$. Here, we are also inspired by the
results about the integrals on finite-dimensional Hopf algebras as
found in [VD3].
\snl
We first need to notice that the definition of a Hopf
group-coalgebra as given by Virelizier in [V, Definition 1.3] is
somewhat more general than the original definition given by Turaev
in [T, Definition 11.2]. However, in the finite-dimensional case
that we consider here, Virelizier showed that the two notions
coincide  (see [V, Corollary 3.7]). His proof of this result makes
use of the integral on the Hopf group-coalgebra. We will give a
more direct proof of this result in Lemma 3.2 below.
\nl
For convenience of the reader, we recall the following definition,
given by Virelizier (Definition 1.3 in [V]).

\inspr{3.1} Definition \rm
Let $G$ be a group. Consider a family $(A_p)_{p\in G}$ of
finite-dimensional unital algebras. Suppose we have the following.
\item{} i) For all $p,q\in G$, we have a unital homomorphism
$\co_{p,q}:A_{pq}\to A_p\ot A_q$ satisfying 'coassociativity' in
the sense that
$$(\co_{p,q}\ot\iota)\co_{pq,r}=(\iota\ot\co_{q,r})\co_{p,qr}$$
on $A_{pqr}$ for all $p,q,r$ in $G$.
\item{} ii) There is a unital homomorphism $\varepsilon_e:A_e\to \Bbb C$
such that for all $p$
$$(\iota\ot\varepsilon)\co_{p,e}(a)=(\varepsilon\ot\iota)\co_{e,p}(a)=a$$
when $a\in A_p$.
\item{} iii) There exists a {\it linear map} $S_p:A_p\to
A_{p^{-1}}$ for all $p$ such that for all $a\in A_e$ we get
$$\align
m(S_{p^{-1}}\ot\iota)(\co_{p^{-1},p}(a))&=\varepsilon_e(a)1\\
m(\iota\ot S_{p^{-1}})(\co_{p,p^{-1}}(a))&=\varepsilon_e(a)1
\endalign$$
for all $p$ (where $m$ denotes multiplication in $A_p$).
\snl
With these conditions, the family $(A_p)_{p\in G}$ is called a
'finite-dimensional' Hopf $G$-coalgebra (in the sense of
Virelizier).
\einspr

Remark the similarity with the notions in Theorem 1.5 in Section
1. Also observe that this notion is close to the definition given
by Turaev, the difference being that here it is not assumed that
the maps $S_p$ are anti-isomorphisms. It is just required that
they are linear maps. It is well-known in Hopf algebra theory that
it is sufficient to assume that the antipode is linear. It follows
that it is an anti-homomorphism. This is also the case here. And
because the components are finite-dimensional, we even get also
bijectivity of the maps $S_p$. Again, this is also expected from
Hopf algebra theory (see [A] and [S]). As mentioned before, we
obtain these result directly, in the following lemma.

\inspr{3.2} Lemma \rm
With the notations and assumptions of Definition 3.1, we get that
$S_p:A_p\to A_{p^{-1}}$ is a anti-isomorphism for all $p\in G$.
\bf \snl Proof: \rm Denote by $A_p'$ the dual space of $A_p$ for
all $p$. Consider the direct sum $A^*=\oplus_{p\in G} A_p'$ of
these dual spaces (and consider the components as subspaces). It
is a subspace of the full dual $A'$ of the direct sum
$A=\oplus_{p\in G} A_p$. We will now make $A^*$ into a Hopf
algebra and show that it has integrals. Then we will use results
from usual Hopf algebra theory to prove the lemma.
\snl
First we make $A^*$ into a unital algebra. We define the product
$fg$ of elements $f\in A_p'$ and $g\in A_q'$ by $(fg)(x)=(f\ot
g)(\co_{p,q}(x))$ for $x\in A_{pq}$ and we put $fg$ equal to $0$
on other components. So in fact, $fg\in A_{pq}'$. The product is
associative because of condition i) in Definition 3.1. The element
$\varepsilon_e$ in $A_e'$ is an identity for $A^*$ because of
condition ii) in Definition 3.1.
\snl
Now, we define a coproduct $\co$ on this algebra $A^*$ by
dualizing the product. We let $\co:A_p'\to A_p'\ot A_p'$ be given
by $(\co(f))(x\ot y)=f(xy)$ when $x,y\in A_p$. This map is a
unital homomorphism and it will be coassociative. A counit
$\varepsilon$ on $A^*$ is given by $\varepsilon(f)=f(1)$ where $1$
is the identity in $A_p$ when $f\in A_p'$. An antipode $S$ on
$A^*$ is given by the formula $(S(f))(x)=f(S_{p^{-1}}(x))$ where
$x\in A_{p^{-1}}$ and $f\in A_p'$. So, $S:A_p'\to A_{p^{-1}}'$ for
all $p$. That it is an antipode follows from condition iii) in
Definition 3.1.
\snl
Hence, we have made $A^*$ into a Hopf algebra. Recall that we only
need the antipode to be linear. It automatically follows that it
is a anti-homomorphism.
\snl
It also follows that $\co\circ S=(S\ot S)\circ\co^{\text{op}}$
where $\co^{\text{op}}$ is the opposite comultiplication, obtained
by applying the flip automorphism (see e.g.\ [S, Proposition
4.0.1]). From this it easily follows that the map $S_p$ is a
anti-homomorphism from $A_p$ to $A_{p^{-1}}$ for all $p$.
\snl
It remains to show that the maps $S_p$ are bijections.
\snl
To prove this, first observe that $(A_e,\co_{e,e})$ is a
finite-dimensional Hopf algebra. This is an easy consequence of
the conditions in Definition 3.1. We know that there must be a
left cointegral $h$ in $A_e$. This is an element $h\in A_e$
satisfying $ah=\varepsilon_e(a)h$ for all $a\in A_e$. We will now
use this element to get an integral on $A^*$. We define
$\varphi(f)=f(h)$ for $f\in A_e$ and $\varphi(f)=0$ on other
components. This defines a linear functional on $A^*$. It is not
hard to show that it is a left integral.
\snl
It is known that the antipode in a Hopf algebra with integrals
must be injective (see e.g.\ [S, Corollary 5.1.7]). In particular,
the maps $S_p:A'_p\to A'_{p^{-1}}$ will be injective. And because
this also is the case for $S_{p^{-1}}:A'_{p^{-1}}\to A'_p$, and
the spaces are finite-dimensional, we must have that $A'_p$ and
$A'_{p^{-1}}$ have the same dimensions and that these maps are
bijections. Then also the adjoint maps $S_p:A_p\to A_{p^{-1}}$
will be bijections.
\snl
This completes the proof of the lemma.
\einspr

It is an immediate consequence of this lemma that the direct sum
$A$ of these algebras is a regular $G$-cograded multiplier Hopf
algebra because the assumptions of Theorem 1.5 are now fulfilled.
\nl
So, in the rest of this section, we consider a regular
$G$-cograded multiplier Hopf algebra $A$ with finite-dimensional
unital components $(A_p)_{p\in G}$.
\snl
We get the following result.

\inspr{3.3} Proposition \rm
Assume that $(A,\co)$ is a regular $G$-cograded multiplier Hopf
algebra with finite-dimensional unital components. Then $A$ is a
multiplier Hopf algebra of discrete type.
\snl\bf Proof: \rm
Recall from Theorem 1.5 that the family of components is provided
with these homomorphims $\co_{p,q}:A_{pq}\to A_p\ot A_q$, just
like in the Definition 3.1 above. Again $(A_e,\co_{e,e})$ will be
a finite-dimensional Hopf algebra. As in the proof of the previous
proposition, let $h$ be a left cointegral in $(A_e,\co_{e,e})$.
Because of the structure of the algebra $A$, being a direct sum of
the algebras $(A_p)_{p\in G}$, we consider $h$ as sitting in $A$
and it will still satisfy $ah=\varepsilon(a)h$, for all $a\in A$.
So, $h$ is also a left cointegral for the multiplier Hopf algebra
$(A,\co)$. Then, by definition, this multiplier Hopf algebra is of
discrete type (see Definition 5.2 in [VD4]).
\einspr

The following is then an immediate consequence (see [VD-Z1]).

\inspr{3.4} Corollary \rm
With the notations and assumptions of Proposisition 3.3 we have
the following.
\item{} i) The multiplier Hopf algebra $A$ has integrals (Theorem
2.10 in [VD-Z]).
\item{} ii) The dual $\hat A$ of $A$ is a (usual) Hopf algebra. In
this case, $\hat A=A^*$ where $A^*=\oplus_{p\in G} A_p'$. The Hopf
algebra structure on $A^*$ is as in the proof of Lemma 3.2.
\einspr

Combining the previous results, we recover the result (on the
existence of integrals) of Virelizier (Section 3 in [V]). Recall
that the uniqueness of the integrals is a general property,
already stated in the previous section.
\nl
The rest of this section is now devoted to obtain an explicit
formula for the left integral on this $G$-cograded multiplier Hopf
algebra with finite-dimensional unital components. The formula we
obtain is of the same nature as the formula for the left integral
on a finite-dimensional Hopf algebra as given in [VD3]; see also
Section 5 in [VD4].
\snl
Remark that the dual space $A'$ of $A$ (the space of all linear
functionals on $A$) is given by $\prod_{p\in G}A_p'$ (as we have
seen before, in Section 2). In general, for a multiplier Hopf
algebra $A$, the dual space $A'$ can not be made into an algebra
by dualizing the product. In this case, we have $\hat A=A^*$ with
$A^*=\oplus_{p\in G} A_p'$ as we saw in Corollary 3.4 and we have
an algebra structure on $\hat A$. We can extend this product.
Indeed, it is possible to define $fg\in A'$ for $f,g\in A'$ as
soon as one of the factors actually sits in the reduced dual $A^*$
(in fact, this can be done for any regular multiplier Hopf
algebra, see [VD1]). The product is still given by $(fg)(a)=(f\ot
g)\co(a)$, an expression that makes sense if e.g.\ $f\in A^*$
because then the first leg of $\co(a)$ will be covered.
\snl
Using this definition, it is easily seen that a linear function
$\varphi$ on $A$ will be left invariant if and only if
$f\varphi=f(1_p)\varphi$ for all $f\in (A_p)'$ and for all $p\in
G$. We will use this to prove left invariance below.
\snl
Recall that all the components $A_p$ are assumed to be
finite-dimensional. Therefore we can introduce the following.

\inspr{3.5} Notations \rm
For each $p\in G$, let $(e_{p,i})_i$ be a basis in $A_p$ and
$(f_{p,i})_i$ a dual basis in the dual space $A_p'$. For $f\in A'$
and $a\in A_p$, denote by $f(\,\cdot\,a)$ the linear functional on
$A_p$ given by $x\mapsto f(xa)$.
\einspr

We need the following lemmas.

\inspr{3.6} Lemma \rm
With the notations and assumptions from before, we have the
following two formulas.
\item{} i) For all $p,q$ we have
$$\sum_i f_{(pq),i}\ot
\co_{p,q}(e_{(pq),i})=\sum_{j,k}f_{p,j}f_{q,k} \ot e_{p,j} \ot
e_{q,k}$$ in $A^*\ot A_p\ot A_q$.
\item{} ii) For all $p$ we have
$$\sum_i \co(f_{p,i})\ot
e_{p,i}=\sum_{j,k}f_{p,j}\ot f_{p,k} \ot e_{p,j} e_{p,k}$$
 in $A^*\ot A^*\ot A_p$.
\einspr

The proof is straightforward. It uses the facts that the product
in $A^*$ is dual to the coproduct in $A$ (to prove i)) and that
the coproduct in $A^*$ is dual to the product in $A$ (to prove
ii)). Of course, also the basic property of dual bases has to be
used several times. It says that for all $p$ and all $a\in A_p$ we
have $\sum_if_{p,i}(a)e_{p,i}=a$.

\inspr{3.7} Lemma \rm
Again with the notations and assumptions from before, we get the
following.
\item{} i) For all $f\in A^*$, $g\in A'$ and $a\in A$ we have
$$(fg)(\,\cdot\,
S(a))=\sum_{(a)}f(\,\cdot\,S(a_{(2)}))g(\,\cdot\,S(a_{(1)})).$$
\item{} ii) For all $f\in A'$, $a\in A$ and $p\in G$ we have
$$\sum_if_{p,i}(f(\,\cdot\,S(a))) \ot e_{p,i} =
\sum_{i,(a)} (f_{p,i}f)(\,\cdot\,S(a_{(1)})) \ot e_{p,i}a_{(2)}.$$

\snl\bf Proof: \rm
First observe that $f$ is a reduced functional on $A$ and
therefore, the element $S(a_{(2)})$ will be covered in the first
formula. Moreover, the product $fg$ is well defined as a product
of an element in $A^*$ with an element in $A'$. The result is in
$A'$, but as we include $S(a)$ in the argument, the functional in
the left hand side of the first formula is an element in $A^*$. In
the right hand side of this formula, we have a sum of products of
elements in $A^*$ and so the result is again in $A^*$. So, the
first formula in the lemma makes sense within $A^*$. The proof of
i) is an immediate consequence of the definition of these
products.
\snl
Next we consider ii). Now, $a_{(2)}$ is covered by the elements
$e_{p,i}$. Similarly as in i), we have products of functionals
with at least one factor in $A^*$ and anyway, in both sides of the
equation, we have a result with functionals in $A^*$. Using i) the
right hand side of the formula in ii) can be written as
$$\sum_{i,(a)} f_{p,i}(\,\cdot\,S(a_{(2)}))(f(\,\cdot\,S(a_{(1)})) \ot e_{p,i}a_{(3)}.$$
Now we use ii) of Lemma 3.6 to see that this expression gives
$$\align
\sum_{i,j,(a)}f_{p,j}(S(a_{(2)}))(f_{p,i}(f(\,\cdot\,S(a_{(1)}))))
  & \ot e_{p,i}e_{p,j}a_{(3)}\\
  & =\sum_{i,(a)}f_{p,i}(f(\,\cdot\, S(a_{(1)}))) \ot
  e_{p,i}S(a_{(2)})a_{(3)} \\
  & = \sum_i f_{p,i}(f(\,\cdot\,S(a)))\ot e_{p,i}.
\endalign$$
This proves ii).
\einspr

Now that we have obtained the main technical results, we will
define a candidate for the left integral. The previous lemmas will
make it rather easy to prove that we do have a non-zero left
invariant functional, defined as follows.

\inspr{3.8} Definition \rm
Let $f$ be any linear functional in $A_e'$. Define a linear
functional $\varphi_f$ on $A$ by $\varphi_f(a)=\sum_{p,i}
(f_{p,i}f)(aS^2(e_{p,i}))$.
\einspr

Observe that, because $f\in A_e'$, the products $f_{p,i}f$ are
again functionals on $A_p$ and that $S^2$ leaves $A_p$ globally
invariant so that $aS^2(e_{p,i})\in A_p$. Because $a\in
\oplus_{p\in G} A_p$, the element $a$ will only have finitely many
components and therefore, we have a sum over only finitely many
$p$ above. So, this definition makes sense.
\snl
Compare this formula with the one in Proposition 1.1 in [VD3]. It
is a dual version of this formula, generalized to the case we
study here.
\snl
Therefore, the following should not come as a surprise.

\inspr{3.9} Theorem \rm
Let $(A,\co)$ be a regular $G$-cograded multiplier Hopf algebra
with finite-dimensional unital components. For any $f\in A_e'$,
the linear functional $\varphi_f$ on $A$, as defined in Definition
3.8, is left invariant on $(A,\co)$.

\snl\bf Proof: \rm
As explained before, to prove left invariance of $\varphi_f$, we
need to show that $g\varphi_f=g(1_p)\varphi_f$ in $A'$ for all
$p\in G$ and $g\in A_p'$. This is equivalent with requiring
$\sum_if_{p,i}\varphi_f \ot e_{p,i}=\varphi_f \ot 1_p$ in $A'\ot
A_p$ for all $p$. We use $1_{p}$ to denote the identity in
$A_{p}$.
\snl
Now, let $p,q\in G$. We have, using the definition of $\varphi_f$,
$$\align \sum_i (f_{p,i}\varphi_f)(\,\cdot\,1_{pq}) \ot e_{p,i}
     & = \sum_i f_{p,i}(\varphi_f(\,\cdot\,1_q)) \ot e_{p,i}\\
     & = \sum_{i,j} f_{p,i}((f_{q,j}f)(\,\cdot\,S^2(e_{q,j}))) \ot
     e_{p,i}.
\endalign$$
By Lemma 3.7.ii), this last expression equals
$$\sum_{i,j,(e_{q,j})} (f_{p,i}f_{q,j}f)(\,\cdot\,S^2({e_{q,j}}_{(2)}))
    \ot e_{p,i}S({e_{q,j}}_{(1)}).$$
Then, using Lemma 3.6.i) we get
$$\align\sum_{i,j,k}(f_{p,i}f_{p^{-1},j}f_{pq,k}f)(\,\cdot\,S^2(e_{pq,k}))
    &\ot e_{p,i}S(e_{p^{-1},j})\\
    &=\sum_k(f_{pq,k}f)(\,\cdot\, S^2(e_{pq,k}))\ot 1_p\\
    &=\varphi_f(\,\cdot\,1_{pq})\ot 1_p.
\endalign$$
This is true for all $q$ and hence we obtain the required formula
$\sum_if_{p,i}\varphi_f \ot e_{p,i}=\varphi_f \ot 1_p$ in $A'\ot
A_p$ for all $p$. This completes the proof.
\einspr

In order to get a left integral on $A$ with this formula, we need
to show that for at least one $f$, the functional $\varphi_f\neq
0$. By the uniqueness of left integrals, then we know that
$\varphi_f$ will be a scalar multiple of the left integral. This
problem is similar as with the formula that we used in [VD3]. We
solve it in the following proposition.

\iinspr{3.10} Proposition \rm
Take $p\in G$ and any $g\in A_p'$. Then
$g(a)=\sum_i\varphi_{f_{p^{-1},i}g}(aS(e_{p^{-1},i}))$ for all
$a\in A_p$. Consequently, $\varphi_h\neq 0$ for some $h\in A_e'$.

\snl\bf Proof: \rm
From equation i) in Lemma 3.6 we get
$$\varepsilon_e\ot 1_p=
\sum_{i,j}(f_{p,j}f_{p^{-1},i}\ot S(e_{p^{-1},i})S^2(e_{p,j}))$$
for all $p$. Therefore, any $g\in A_p'$ has the form
$$g=\sum_{i,j}(f_{p,j}f_{p^{-1},i}g)(\,\cdot\,S(e_{p^{-1},i})
S^2(e_{p,j}))$$ or equivalently
$$g=\sum_i\varphi_{f_{p^{-1},i}g}(\,\cdot\,S(e_{p^{-1},i})).$$
Now, it follows from this formula that, if $g\neq 0$, we must have
that for some $i$ also $\varphi_h\neq 0$ where $h=f_{p^{-1},i}g$.
\einspr

So, we have shown that $\varphi_f$ will be non-zero for at least
some $f\in A'_e$. We proved already that any $G$-cograded
multiplier Hopf algebra with finite-dimensional unital components
has integrals (Proposition 3.3 and Corollary 3.4). Now, we also
have a constructive proof of this result.

\nl\nl

\bf References \rm
\bigskip
\parindent 1 cm
\item{[A]} E.\ Abe: \it Hopf algebras. \rm Cambridge University Press (1977).
\smallskip
\item{[D]} L.\ Delvaux: {\it Semi-direct products of multiplier
Hopf algebras: smash coproducts.} Commun.\ Algebra 30 (2002),
5979--5997.
\smallskip
\item{[D-VD]} L.\ Delvaux \& A.\ Van Daele: {\it The Drinfel'd double
for goup-cograded multiplier Hopf algebras.} Preprint L.U.C.\ and
K.U.Leuven (2004). Server version math.QA/04????.
\smallskip
\item{[H-A-M]} A.S.\ Hegazi, A.T.\ Abd El-Hafez \& M.\ Mansour:
{Multiplier Hopf group coalgebras}. Preprint Mansoura University
(2002).
\smallskip
\item{[K]} J.\ Kustermans: {\it The analytic structure of
algebraic quantum groups.} J. Algebra 259 (2003), 415-450.
\smallskip
\item{[L-VD]} M.B.\ Landstad \& A.\ Van Daele: {\it (In
preparation)}. Preprint NTNU Trondheim and K.U.Leuven (2004).
\smallskip
\item{[S]} M.E.\ Sweedler: \it Hopf algebras. \rm Mathematical Lecture Note
Series. Benjamin (1969).
\smallskip
\item{[T]} V.G.\ Turaev: {\it Homotopy field theory in dimension 3 and crossed
group-categories.} IRMA Preprint, Strasbourg, 2000. Server version
math.GT/0005291.
\smallskip
\item{[VD1]} A.\ Van Daele: {\it Multiplier Hopf algebras.} Trans.\
Amer.\ Math.\ Soc.\ 342 (1994), 917-932.
\smallskip
\item{[VD2]} A.\ Van Daele: {\it Discrete quantum groups.} J.\ of Alg.\
180 (1996), 431-444.
\smallskip
\item {[VD3]} A.\ Van Daele: {\it The Haar measure on finite quantum
groups}. Proc.\ Amer.\ Math.\ Soc.\ 125 (1997), 3489-3500.
\smallskip
\item{[VD4]} A.\ Van Daele: {\it An algebraic framework for group duality.}
Adv. in Math.  140 (1998), 323--366.
\smallskip
\item{[VD-Z1]} A.\ Van Daele \& Y.\ Zhang : {\it Multiplier Hopf
algebras of discrete type}. J.\ of Alg.\ 214 (1999), 400-417.
\smallskip
\item{[VD-Z2]} A.\ Van Daele \& Y.\ Zhang:
{\it A survey on multiplier Hopf algebras.} Proceedings of the
conference in Brussels on Hopf algebras. Hopf Algebras and Quantum
Groups, eds. Caenepeel/Van Oystaeyen (2000), 269-309. Marcel
Dekker (New York).
\smallskip
\item{[V]} A.\ Virelizier: {\it Hopf group-coalgebras.} IRMA
Preprint, Strasbourg, 2000. Server version math.QA/0012073. J.\
Pure Appl.\ Algebra 171 (2002),75--122.
\smallskip
\item{[W1]} S.\ Wang: {\it Group twisted smash products and
Doi-Hopf modules for T-coalgebras.} Preprint (2003). To appear in
Communications in Algebra.
\smallskip
\item{[W2]} S.\ Wang: {\it Group entwining structures and
group-coalgebra Galois extensions.} Pre\-print (2004). To appear
in Communications in Algebra
\smallskip
\item{[Z]} M.\ Zunino: {\it Double construction for crossed Hopf
coalgebras.} IRMA Preprint, Strasbourg, 2002. Server version
math.QA/0212192.
\smallskip

\end

\end